\documentclass[pdflatex,sn-mathphys]{sn-jnl}

\jyear{2022}%

\theoremstyle{thmstyleone}%
\theoremstyle{thmstyletwo}%

\theoremstyle{thmstylethree}%

\usepackage{fancyvrb}
\usepackage{relsize}
\usepackage{setspace}

\begin{document}

\title{Who Finds the Short Proof?}
\subtitle{An Exploration of Boolos' Curious Inference 
using Higher-order Automated Theorem Proving}



\author*[1,2]{\fnm{Christoph} \sur{Benzmüller}}\email{christoph.benzmueller@uni-bamberg.de}

\author[1,2]{\fnm{David} \sur{Fuenmayor}}\email{david.fuenmayor@uni-bamberg.de}

\author[3]{\fnm{\\Alexander} \sur{Steen}}\email{alexander.steen@uni-greifswald.de}

\author[4]{\fnm{Geoff} \sur{Sutcliffe}}\email{geoff@cs.miami.edu}

\affil[1]{\orgdiv{AI Systems Engineering}, \orgname{Otto-Friedrich-University Bamberg},\\ \orgaddress{\street{An der Weberei 5}, \city{Bamberg}, \postcode{96049}, 
\country{Germany}}} 

\affil[2]{\orgdiv{Department~of Mathematics and Computer Science}, \orgname{FU Berlin},\\ \orgaddress{\street{Arnimallee 7}, \city{Berlin}, \postcode{14195}, 
\country{Germany}}}

\affil[3]{\orgdiv{Institute of Mathematics and Computer Science}, \orgname{University of Greifswald}, \orgaddress{\street{Walter-Rathenau-Str. 47}, \city{Greifswald}, \postcode{17489}, 
\country{Germany}}}

\affil[4]{\orgdiv{Department of Computer Science}, \orgname{University of Miami},\\ \orgaddress{\street{1365 Memorial Drive}, \city{Coral Gables}, \state{FL} \postcode{33124-4245}, \country{USA}}}


\abstract{This paper reports on an exploration of Boolos' Curious Inference, using higher-order automated theorem provers (ATPs). Surprisingly, only suitable shorthand notations had to be provided by hand for ATPs to find a short proof. The higher-order lemmas required for constructing a short proof are automatically discovered by the ATPs. Given the observations and suggestions in this paper, full proof automation of Boolos' and related examples now seems to be within reach of higher-order ATPs.}




\maketitle

\section{Which Automated Theorem Prover to choose for 
Really Difficult Problems?}
\label{section1}

Consider the following thought experiment: Shortly after their death, Folbert and Holly arrive at the gates of heaven, and both want to enter. Unfortunately only one of the two can be admitted, and they have  to settle this in a little contest: Folbert and Holly are both asked to (i)~choose either a first-order (FO) automated theorem prover (ATP), or, alternatively, an higher-order (HO) ATP; and then (ii)~pose a FO proof problem, encode it in FO logic, and give it to the other's ATP. 
The one whose ATP solves the proof problem posed to it the faster will be admitted to heaven. 
There is a time limit for the contest -- until midnight the same day. If the battle ends in a draw because neither ATP finds a proof by midnight, the contest is repeated the next day, and so on. 
Which ATP should Folbert and Holly choose, and which proof problem? 
Is there a winning strategy?

FOlbert turns out to be a FO logic enthusiast, and announces that he will use a FO ATP. HOlly likes HO logic better, and chooses an HO ATP. Which exact system they choose is not revealed to the other person, nor how it functions internally. 
At first glance, Folbert seems to be in a more promising position because the proof problems they must provide have to be encoded in FO logic. 
However, as will be seen, only Holly has a realistic chance to win this contest (maybe not on the first day though), provided she chooses a suitable proof problem. 

Key to Holly's advantage are the (hyper-)exponentially shorter proofs that are possible as one moves up the ladder of expressiveness from first-order logic to second-order logic, to third-order logic, and so on \cite{GoedelProofLength}. 
The fact that the proof problems are stated in FO logic does not matter. 
When stating the same problem in the same FO way but in higher-order logic, much shorter proofs are possible, some of which might be (hyper-)exponentially shorter than the proofs that can be found with FO ATPs. 
A very prominent example of such a short proof is that of \textit{Boolos' Curious Inference} \cite{BCI}. 

The key to the shorter proofs in HO logic is the possibility of introducing powerful shorthand notations in which structural aspects of the problem at hand can be expressed in ways that are not possible in FO logic. 
The use of HO variables and quantifiers in these shorthand notations provides an opportunity to shorten proof arguments by explicitly talking about and 
working with the mathematical structures involved. 
In Boolos' Curious Inference this concerns an inductively defined, very fast growing,  Ackermann-style function. 
A fascinating perspective on meta-level aspects of an otherwise mathematically boring, repetitive FO inference can be provided and exploited this way.
A most interesting aspect of this paper is that the required ``curious'' lemmas, supplied by hand by the ingenious Boolos, can now be synthesized by modern HO ATPs.

Based on our experiments we argue that there is a winning strategy in the contest for only Holly. 
Holly is the only person who is able to reliably prevent her opponent from winning, provided she chooses Boolos' challenge problem to be given to Folbert. 
Folbert can of course try something similar, but this cannot prevent Holly's victory in the long run, since full proof automation of Boolos' challenge problem is now within reach of HO ATPs. 
As the experiments reported in this paper show, significant progress has been achieved in HO ATP systems in recent years, which enables the automatic exploration of variants of Boolos' Curious Inference.
In particular, powerful and intuitive lemmas can be automatically discovered by HO ATPs, and used in short proofs. 
The only thing left for the human to do has been to suggest suitable shorthand notations.

In Section~\ref{sec:ProofByE} the proofs produced by the E ATP \cite{Eprover}, which 
are based on three lemmas, are discussed in some detail (we chose E because its proofs were the most readable and because E also performed best). 
Further proofs are reported by the HO ATPs 
cvc5 \cite{CVC5},
Ehoh \cite{Ehoh},
Leo-III \cite{Leo3}, and
Zipperposition \cite{Zipperpin}; these contributions will be studied in more depth in future work.
This paper also discusses the remaining key challenge that needs to be addressed by HO ATPs to enable full automation of finding a solution to Boolos' and similar problems: constructing the required shorthand notations using controlled cut/comprehension introduction, as was already hinted at in \cite{LostProof}.

 Paper outline: Section~\ref{BCI} briefly recaps Boolos' challenge problem, and points to some formalisations that have been verified with interactive HO proof assistants. 
 Section~\ref{sec:ProofByE} demonstrates that HO ATPs can easily solve Boolos' problem and produce short proofs when suitable shorthand notations are provided.
 Section~\ref{sec:ProofByE} presents and discusses the proofs found by E. 
 Section~\ref{sec:Survey} provides a survey of the results from several HO ATPs. 
 Section~\ref{sec:LostProof} explains why the required shorthand notation cannot (yet) be generated with state-of-the art HO ATPs, and points to interesting future research on controlled cut/comprehension introduction.
 Section~\ref{Conclusion} concludes the paper: Holly is the only one who will one day enter heaven. 

\section{Boolos' Curious Inference} 
\label{BCI}

In his article ``A Curious Inference'' \cite{BCI}, Boolos presents the following proof problem, consisting of axioms A1-A5 with the conjecture C. This challenge is referred to as BCP (for Boolos' Curious Problem) in the rest of this paper: 
\begin{align}
&   \forall n. f(n,e) = s(e)  \tag{A1} \\
&   \forall y. f(e,s(y)) = s (s (f(e,y))) \tag{A2}  \\
&    \forall x\,y. f(s(x),s(y)) = f(x,f(s(x),y))  \tag{A3} \\
&    d(e)  \tag{A4} \\
&   \forall x. d(x) \rightarrow  d(s(x)) \tag{A5}  \\[.5em]
&   d(f(s(s(s(s(e)))),s(s(s(s(e)))))) \tag{C} 
\end{align}
In the first three axioms there are three uninterpreted FO constant symbols: a nullary FO constant symbol $e$ (intuitively, think of it as the number \textit{one}, if it helps), an unary FO function symbol $s$ (think of is as the \textit{successor} function), and a binary FO function symbol $f$, whose semantics is inductively characterised in the axioms A1-A3 (where A2 and A3 are recursive; note that A3 is actually double recursive). 
These axioms capture the fact that $f$ belongs to a class of extremely fast growing functions, also known as Ackermann(-style) functions~\cite{ackermann1928}. 
Exhaustive evaluation of the term $f(s(s(s(s(e)))),s(s(s(s(e)))))$ with these recursive equations unfolds it to a term that contains more $s$'s than there are atoms in the universe. 
The conjecture C formulated by Boolos is to prove that $d$ (think of it as the \textit{is-a-natural-number} predicate) holds for $f(s(s(s(s(e)))),s(s(s(s(e)))))$. $d$ is actually an uninterpreted FO predicate symbol that holds for $e$ (A4), and is propagated from any $x$ to $sx$ (A5). A formalisation of BCP in the TPTP THF language \cite{Sut22-IGPL,THF}
can be found in Appendix \ref{sec:Boolos1.p}.

It is easy to see that this proof problem is solvable in FO logic by applying an astronomically large number of modus ponens steps to A4 and (instances of) A5. 
Boolos showed that for a cut-free FO calculus the magnitude of this number is comparable to $f(5,5)$ for the Ackermann function $f$ interpreted on the naturals. 
Due to the enormous number of required inference steps it is certain that no such proof can realistically be found, let alone represented, within cut-free FO proof systems. 
However, Boolos also showed that there exists a short proof of BCP in second-order logic. 
This proof has been encoded and verified in the mathematical proof assistant systems Omega \cite{Omega} and Mizar \cite{Mizar} by Benzmüller and Brown  \cite{BoolosOmegaMizar}, at essentially the same level of granularity as that of Boolos' proof. 
This earlier work was recently repeated by Ketland \cite{KetlandAFP} in the interactive proof assistant Isabelle/HOL \cite{Isabelle}.

The short proof by Boolos makes use of specific instances of the comprehension principle of second-order logic, an axiom schema given by
\begin{align}
& \exists R. \forall x_1 \ldots x_n.\; R(x_1, \ldots, x_n) \leftrightarrow  \varphi(x_1, \ldots, x_n)   \tag{COM}
\end{align}
where $\varphi(x_1, \ldots, x_n)$ is a second-order formula with $x_1, \ldots, x_n$ among its free variables, 
and $R$ is a second-order variable not free in $\varphi$.
The comprehension principle postulates the existence of relations (or predicates) that can be defined within the second-order language, or, equivalently, that there exists some $R$ that can act as a shorthand notation for the property expressed by $\varphi(x_1, \ldots, x_n)$. 
For function symbols a similar axiom scheme exists. 
The idea is to choose the right instances of COM  such that helpful lemmas can be used to enable a short (second-order) proof for BCP. Details are provided in Boolos' paper. 

\sloppy In a HO logic based on  $\lambda$-terms, such as Church's type theory \cite{Church40,Andrews,SEP}, the existence of $R$ in COM is generally guaranteed: simply choose $R = \lambda x_1, \ldots, x_n.\,\varphi(x_1, \ldots, x_n)$. 
The COM instances are easily provable, and axiomatisations of the COM principle are generally avoided in HO ATPs. 
During HO proof search, $\lambda$-terms for $R$ are often synthesised on the fly by HO unification. 
In less trivial cases, however, parts of the term $R$ need to be ``guessed'' by the ATPs, e.g., by applying a technique called \textit{primitive substitution}. 
A well known case is Cantor's theorem, where an initial part of the diagonal set description needs to be guessed before the rest can be synthesised; see, e.g.,~\cite{Andrews89,Leo2} for further details on primitive substitution.

\section{Proof Variants} 
\label{sec:ProofByE}

Without further support BCP cannot be solved by today's HO ATP systems.  
Quite surprisingly (to the authors, at least), with a little help provided in 
the form of suitable shorthand notations, the ATPs can automatically discover 
lemmas that enable a short proof.

In the remainder of this paper we use Boolos' convention to write $sssse$ instead of $s(s(s(s(e))))$; analogously for $ssse$, $sse$, $se$, $sx$, etc. Moreover, we generally assume HO notation, so that e.g.~$f(sssse,sssse)$ is represented as $((f\,sssse)\,sssse)$ or simply as $fssssesssse$, avoiding unnecessary parentheses and even spaces when the correct term and formula structures can be easily inferred in context. 

\subsection{Variants using Two Shorthand Notations}
\label{TwoNotations}
\newcommand{\cl}{\texttt{C}}
\newcommand{\githubbaseurl}
{https://github.com/cbenzmueller/BoolosCuriousInference-ATP/tree/main/Boolos1.proof}
\newcommand{\clref}[2]{\href{\githubbaseurl\##2}{\cl_{#1}}\colon}
\newcommand{\clname}[2]{\ensuremath{\href{\githubbaseurl\##2}{\cl_{#1}}}}

When given the following two shorthand notations, various HO ATPs can find a short proof for BCP. 
The HO representation of BCP is augmented with two additional uninterpreted symbols
$ind$ and $p$. 
They are governed by two additional axioms, (Def\_$ind$) and (Def\_$p$), stating equalities between constant symbols and $\lambda$-terms.
%
\begin{align}
& ind = \lambda X.\, Xe \wedge \forall x.\, Xx \rightarrow X sx \tag{Def\_$ind$} \\
& p = \lambda x\,y.\, (\lambda z.\, \forall X.\ ind\ X \rightarrow X z)\ fxy \tag{Def\_$p$} 
\end{align}
Recalling the discussion of COM in Section~\ref{BCI}, note that $ind$ and $p$ are just shorthand notations for the HO $\lambda$-terms given on the right. 
The $\lambda$-term abbreviated by $ind$ expresses what it means to be an inductive set defined over $e$ and $s$. 
The $\lambda$-term abbreviated by $p$ expresses that for any $x$ and $y$ given as arguments, $fxy$ is in the smallest inductive set defined over $e$ and $s$.\footnote{%
It is acknowledged that the use of terminology might be slightly misleading here: $s$ might not be injective, and $e$ might have predecessors; in particular $s$-cycles have not been axiomatised away. 
There is thus an alternative explanation and terminology that could be used: 
The declarations of the constants $e::i$ and $s::i\rightarrow i$ express (semantically speaking) that every model behaves as an algebra over the signature $\Sigma = \{e/0$, $s/1\}$. Then $ind~X$ means that $X$ is a $\Sigma$-(sub)algebra of the model. 
Thus, the term ``smallest inductive set defined over $e$ and $s$'' (suggestively called ``$N$'' by Boolos) corresponds to the infimum/intersection of all $\Sigma$-(sub)algebras. 
From universal algebra it is known that the set of all subalgebras of any algebra forms a complete lattice, where infimum corresponds to set-intersection, and thus $N :=\bigcap~ind$ is guaranteed to exist. Hence $p$ can be alternatively written as $\lambda x\,y.\,N\,(fxy)$. \label{footnote}}
By adding these shorthand notations HO ATPs can find short proofs for BCP. 
In Section~\ref{sec:Survey} the experiments conducted using the TPTP THF infrastructure for HO logic \cite{TPTP,THF} are summarized. 
Interesting lemmas are automatically discovered, proven, and used. 
No mathematical or logical ingenuity is needed. 
The proof found by E can easily be converted into the following quite intuitive informal proof, divided into four parts. 
The first three parts \textit{I--III} introduce and prove the key lemmas. This was done fully automatically by E. 
In the last part \textit{IV}, the final refutation argument is automatically constructed. 
The new lemmas are used in this proof to derive a contradiction from the negated proof statement. 
The clause names mentioned in the text below reference the annotated formula names in the proofs shown in Appendices \ref{sec:Boolos1IDV.pdf} and \ref{sec:Boolos1.proof}, where the proof graphs and the proofs generated by E are presented.\footnote{%
In order to improve readability, $\cl_i$ is used as a label for the clause named \texttt{c\_0\_$i$} in the formal proof by E. 
Labels are prepended to clauses with a colon, as in $\cl_i \colon \varphi$, where $\varphi$ is some formula. 
Each clause reference is a hyperlink pointing to the clause representation in the original proof output of E (in a Github repository).}

\paragraph{Part I}
 From A4, A5, and Def\_$ind$, it follows that $d$ is an inductive set:
\begin{align*}
& \clref{48}{L381-L384} ind\,d  \tag{L1a}
\end{align*}
From A1, Def\_$ind$, and Def\_$p$, it follows that: 
\begin{align*}
& \clref{31}{L276-L278} pxe \text{ for all $x$} \tag{L1b}
\end{align*}
From A1, A2, Def\_$ind$, and Def\_$p$, it follows that $p e$ is an inductive set:
\begin{align*}
& \clref{59}{L442-L445} ind\, pe \tag{L1c}
\end{align*}
(Remember that p is a binary predicate. Modulo currying and the identification of sets with their characteristic functions, this means that $pe$, and its $\eta$-expansion $\lambda z.\, p e z$, denotes an unary predicate on entities.)

\textit{Proof.} 
Fully formal proofs are provided in Appendix \ref{sec:Boolos1.proof}.

The proof of L1a (see the derivation of clause \clname{48}{L381-L384} from 
A4, A5, and Def\_$ind$) is obvious. 

Proof of L1b (see the derivation of clause \clname{31}{L276-L278} from 
axiom A1, Def\_$ind$, and Def\_$p$): 
From Def\_$ind$ get $\clref{23}{L227-L231} ind\,U \rightarrow U e$  and $\clref{18}{L197-L202} ind\,U \rightarrow  (U x \rightarrow U sx)$, for all $U$ and $x$. Moreover, from Def\_$p$ infer that for all $x$ and $y$ there exists a property $qxy$ (depending on $x$ and $y$) encoded by \texttt{epred1\_2 @ X @ Y} in Appendix \ref{sec:Boolos1.proof}, such that $\clref{16}{L186-L190} qxy\,(fxy) \rightarrow pxy$ and $\clref{19}{L204-L208} ind\,qxy \vee pxy$. 
From \clname{19}{L204-L208} and \clname{23}{L227-L231} it follows $\clref{25}{L238-L242} qxy\,e \vee pxy$, for all $x$ and $y$. From 
\clname{19}{L204-L208} and \clname{18}{L197-L202}, for all $x$, $y$ and $z$, get $\clref{21}{L215-L220} (qxy\,z \rightarrow qxy\,sz) \vee pxy$.
From A1 and \clname{16}{L186-L190} get $\clref{28}{L258-L262} qze\,se \rightarrow pze$, for all $z$. Finally obtain L1b, i.e., clause \clname{31}{L276-L278}, from \clname{28}{L258-L262}, \clname{25}{L238-L242}, and \clname{21}{L215-L220}.

Proof of L1c (see the derivation of clause \clname{59}{L442-L445} from 
A1, A2, Def\_$ind$, and Def\_$p$): 
Instantiating $z$ with $s(fez')$ in \clname{21}{L215-L220} and applying the equation $fesz' = ssfez'$, which follows from A2, get $\clref{50}{L393-L398} (qxy(sfez)' \rightarrow qxy(fesz')) \vee pxy$, for all $z'$. 
From this and \clname{16}{L186-L190} obtain $\clref{53}{L413-L417} qesz'(sfez') \rightarrow pesz'$ and also $\clref{55}{L423-L427} qesz' (fez') \rightarrow pesz'$, both for all $z'$.
Now, from A1 and \clname{16}{L186-L190} get $\clref{28}{L258-L262} qze\, se \rightarrow pze$, for all $z$. 
Moreover, Def\_$ind$ implies that for all properties $U$ there exists some property $kU$ (depending on $U$) encoded by \texttt{esk1\_1 @ U} in Appendix \ref{sec:Boolos1.proof}, such that $\clref{27}{L251-L256} Ue \rightarrow (ind\,U \vee U kU)$ and $\clref{37}{L312-L317} (U e \wedge U skU) \rightarrow ind\,U$. 
From Def\_$p$ it follows $\clref{26}{L244-L249} (ind\,U \wedge pxy) \rightarrow U fxy$, for all $x$, $y$  and $U$. Using \clname{26}{L244-L249} and \clname{27}{L251-L256}
get that for all  $x$ there exists $(k\,px)$ such that for all $U$ holds $\clref{30}{L268-L274} (ind\,U\wedge pxe) \rightarrow (ind\,px \vee U fx(k\,px))$. 
Apply L1b to infer that for all $x$ there exists $(k\,px)$ such that for all $U$ holds $\clref{34}{L292-L297} ind\,U \rightarrow (ind\,px \vee U fx(k\,px))$. 
Using \clname{55}{L423-L427} it then follows from \clname{34}{L292-L297} that $\clref{57}{L433-L436} ind\, pe \vee pes(k\,pe)$.
From \clname{57}{L433-L436} and \clname{37}{L312-L317} finally obtain L1c, i.e., clause \clname{59}{L442-L445}. \qed

\paragraph{Part II}
From A1, A2, A3, Def\_$ind$, and Def\_$p$, it follows: 
\begin{align*}
& \clref{47}{L373-L379} \forall x. \forall Y.\, (ind\,px \wedge ind\, psx \wedge ind\,Y) \rightarrow Y\,fsxsssse \tag{L2}
\end{align*}
Informally, for all entities $x$ and sets $Y$, if $px$, $psx$, and $Y$ are inductive sets (over e and s), then $f\,sx\,sssse$ is in $Y$.

\textit{Proof.} See the derivation of clause \clname{47}{L373-L379} from clauses \clname{8}{L111-L118}--\clname{11}{L132-L140}, corresponding to Def\_$p$, A2, A1, and Def\_$ind$, in Appendix \ref{sec:Boolos1.proof}: 
From A1 and A2 infer $\clref{17}{L192-L195} fese = ssse$, so that with \clname{16}{L186-L190} it follows $\clref{20}{L210-L213} qese\,ssse \rightarrow pese$, where $q$ is again encoded by \texttt{epred1\_2} in Appendix C.
From this and \clname{21}{L215-L220} obtain $\clref{22}{L222-L225} qese\,sse \rightarrow pese$ and $\clref{24}{L233-L236} qese\,se \rightarrow pese$, which together with \clname{25}{L238-L242} leads to $\clref{29}{L264-L266} pese$. 
Using \clname{26}{L244-L249}, i.e., Def\_$p$, and $fese = ssse$ obtain $\clref{32}{L280-L284} ind\, U \rightarrow u(ssse)$, for all $U$. 
Further application of \clname{26}{L244-L249} leads to
$\clref{36}{L305-L310} (ind\,U \wedge ind\,px) \rightarrow U fxssse$ and subsequently to $\clref{43}{L347-L353} (ind\,U \wedge ind\,px \wedge ind\,py) \rightarrow U fxfyssse$. 
Finally, since $fsxsy = fxfsxy$ by A3, obtain L2, i.e., clause \clname{47}{L373-L379}.
\qed

\paragraph{Part III}
From A1, A3, Def\_$ind$, and Def\_$p$, it follows: 
\begin{align*}
& \clref{52}{L406-L411} \forall x.\,~ind\,px \rightarrow ind\,psx  \tag{L3} 
\end{align*}
Informally, for any $x$: if $px$ is an inductive set (over $e$ and $s$), then so is $psx$.

\textit{Proof.}
See the derivation of clause \clname{52}{L406-L411} from clauses \clname{8}{L111-L118}, \clname{10}{L126-L130}, \clname{11}{L132-L140}, and \clname{35}{L299-L303}, corresponding to Def\_$p$, A1, Def\_$ind$, and A3, in Appendix \ref{sec:Boolos1.proof}:
From A1, Def\_$ind$, and Def\_$p$ lemma L1b has been established, given by $\clref{31}{L276-L278} pxe$, for all $x$.
From Def\_$ind$ it follows that $\clref{37}{L312-L317} (U e \wedge U skU) \rightarrow ind\,U$ for all properties $U$.
Clauses \clname{26}{L244-L249} and \clname{34}{L292-L297} have already been inferred above. 
Together with A3, it follows that $\clref{45}{L361-L367} (ind\,px \wedge ind\,U) \rightarrow (U f sx\, s(k\,psx) \vee ind\,psx)$ for every $x$ and every property $U$.
From clause \clname{45}{L361-L367} and Def\_$p$ it follows that $\clref{49}{L386-L391} ind\,px \rightarrow (p\, sx\, s(k\,psx) \vee ind\,psx)$ for all $x$, where $k$ is again the Skolem function called \texttt{esk1\_1} in Appendix \ref{sec:Boolos1.proof}.
An application of \clname{49}{L386-L391} with \clname{37}{L312-L317} yields $(ind\,px \wedge psxe) \rightarrow ind\,psx$, i.e., choose $U = psx$. With a simple application of lemma L1b, instantiated for $sx$, obtain L3, i.e., clause \clname{52}{L406-L411}.
\qed

\paragraph{Part IV: Proof of C by reductio ad absurdum}
To prove C assume $\neg$C and derive a contradiction. 
Assume $\neg d\,fssssesssse$ (clauses \clname{42}{L343-L345} and \clname{46}{L369-L371}). 
From this and L1 obtain that $psssse$ is not an inductive set or that $pssse$ is not an inductive set, i.e., clause \clname{51}{L400-L404}. 
Use L3 to show that $pssse$, $psse$ and $pse$ are not inductive sets (clauses \clname{54}{L419-L421}, \clname{56}{L429-L431}, and \clname{58}{L438-L440}). 
From these, together with L1 and L3, derive a contradiction. 
By reductio ad absurdum $C$ holds.
\qed 
 
\paragraph{Part IV$'$: Constructive proof of C}
The refutation argument IV can be converted into constructive argument IV$'$, avoiding reduction ad absurdum:
From L1 and L3 $pse$ is an inductive set. 
By further applications of L3, $psse$, $pssse$ and $psssse$ are inductive sets. 
Use L1a and L2 to conclude $d\,fssssesssse$. \qed

\paragraph{Note on shorter proofs}
By systematic experimentation using HO ATPs it is possible to further simplify the proof found by E, and to reduce dependencies. For example (as pointed out by an anonymous reviewer), the $ind\,px$ in lemma $L2$ can be avoided and L2 can be proved from just Def\_$p$ and Def\_$ind$.

\subsection{Variants using One Shorthand Notation}
\label{OneNotation}

It turns out that the two shorthand notations, $ind$ and $p$ as introduced above, are dispensable. $ind$ can be avoided when unfolded (but not $\beta$-reduced) in $p$. 
The HO ATPs solve this alternative problem formulation even quicker than before. Instead of $ind$ and $p$ use:
\begin{align}
& p' = \lambda x\,y.\, (\lambda z.\, \forall X.\ (\lambda Q.\, Qe \wedge \forall w.\, Qw \rightarrow Qsw)\, X \rightarrow X z)\,fxy \tag{Def\_$p'$} 
\end{align}
Alternatively, just the following lemma can be suggested, proven, and subsequently used:
\begin{align}
& \exists p'. p' = \lambda x\,y.\, (\lambda z.\, \forall X.\ (\lambda Q.\, Qe \wedge \forall w.\, Qw \rightarrow Qsw)\, X \rightarrow X z)\,fxy \tag{L}
\end{align}
This corresponds to (and is equivalent to) the following instance L$'$ of COM:
\begin{align}
& \exists p'. \forall x\, y.\,  p'xy \leftrightarrow (\lambda z.\, \forall X.\ (\lambda Q.\, Qe \wedge \forall w.\, Qw \rightarrow Qsw)\, X \rightarrow X z)\,fxy \tag{L$'$}
\end{align}

L and L$'$ illustrate the relationship to comprehension and cut-introduction. 
Once introduced, these lemmas can be easily proven by HO ATPs using HO unification, which simply instantiates the existentially quantified $p'$ with the $\lambda$-term: $\lambda x\,y.\, (\lambda z.\, \forall X.\ (\lambda Q.\, Qe \wedge \forall w.\, Qw \rightarrow Qsw)\, X \rightarrow X z)\,fxy$.
Solving BCP fully automatically using HO ATPs thus boils down to speculating L (or L$'$), proving it, and then proving BCP using L (or L$'$). 
The HO ATPs can use the shorthand notations to automatically discover the required cut-lemmas L1a, L1b, L2, and L3.
The encoding of BCP with only L as an axiom is presented in Appendix \ref{sec:BoolosComp.p}. 
A refutation proof for this formulation of BCP, found by E in a few milliseconds, is shown in Appendix \ref{sec:BoolosComp.proof}.

\section{Results of Experiments with HO ATP's}
\label{sec:Survey}

To assess the performance and robustness of HO ATPs in finding short proofs for BCP, 
experiments with different problem encodings were done, using the ATPs
cvc5 1.0 \cite{CVC5},
E 3.0 \cite{Eprover},
Leo-III 1.7.0 \cite{Leo3},
Vampire 4.7 \cite{KV13},
and
Zipperposition 2.1 \cite{Zipperpin}.
These systems were deployed on the StarExec \cite{SST14} Miami cluster running octa-core Intel(R) Xeon(R) E5-2620 v4 @ 2.10GHz CPUs, 128GB memory, and CentOS Linux 7.4.1708 (Core) operating system. A CPU time limit of 300s was set.

The range of CPU times reported below is indicative of the range of difficulty of the problems for ATP.
Comparison of the proof inference statistics between different ATPs is not meaningful, but for the E ATP the statistics indicate that the proofs are of comparable complexity.

\subsection{Using Two Shorthand Notations}

When using the two shorthand notations $ind$ and $p$ discussed in Section~\ref{TwoNotations}, several HO ATPs are able to prove BCP encoding 
given in Appendix \ref{sec:Boolos1.p}. 
It is solved by 
E (13.8s CPU time, 45 inferences in the refutation, 16 inferences deep),
Leo-III (231.0s CPU time, 445 inferences in the refutation, 36 inferences deep),
and
Zipperposition (85.7s CPU time, 29 inferences in the refutation, 6 inferences deep).\footnote{An anonymous reviewer was able to find a proof also with Vampire 4.6, when using a specific parameter setting.}
The proof by E, discussed in Section~\ref{sec:ProofByE}, is presented in Appendix \ref{sec:Boolos1.proof}.

When encoding $ind$ and $p$ as comprehension instances (see Appendix \ref{sec:Boolos1Alt.p}) BCP is proven by only 
Leo-III (215.8s CPU time, 409 inferences in the refutation, 37 inferences deep). 
The resolution proof generated by Leo-III is presented in Appendix~\ref{sec:Boolos1Alt.proof}.

\subsection{Using One Shorthand Notation}

A problem encoding using the comprehension instance $L$ from Section~\ref{OneNotation} as an axiom is shown in Appendix \ref{sec:BoolosComp.p}. 
This problem is solved by only 
E (0.2s CPU time, 53 inferences in the refutation, 20 inferences deep). 
The resolution proof is shown in  Appendix \ref{sec:BoolosComp.proof}.
The related encoding using Def\_$p'$ as an axiom was also proven by only 
E (0.3s CPU time, 47 inferences in the refutation, 17 inferences deep). 

\section{Future Challenge: Shorthand Invention}
\label{sec:LostProof}

Why are the shorthand notations $ind$ and $p$ not found automatically by the HO ATPs? 
Or, alternatively, why is lemma L (or L$'$) not automatically introduced, proven, and then used to obtain a short proofs for BCP, by the HO ATPs?
The answer is that the HO ATPs, in the tradition of FO ATPs, put a strong focus on cut-elimination. 
They do not incorporate  suitable support for controlled cut/comprehension introduction in their proof search. 
While cut-freedom is clearly a desirable property of calculi (a kind of quality seal), things become really interesting, from a mathematical, cognitive, and AI perspective, when cut-elimination is given up at least partially, and some forms of controlled cut-introduction are applied. 
Once controlled cut-introduction steps are applied, powerful further key lemmas can be synthesised automatically within a cut-free resolution-style calculus.
In this sense HO ATPs somehow still neglect a crucial expressivity advantage that they enjoy over FO ATPs.
The situation has already been discussed in a prior paper \cite{LostProof}, but so far not much attention has been paid to this aspect in the HO ATP community. 

The experiments with HO ATPs conducted for this paper, and the discussion above, illustrate the following:
\begin{enumerate}
\item \sloppy Cut-elimination, which has attracted the interest of many ATP researchers and theoreticians, is important for achieving robust proof automation. 
The price of cut-elimination is, however, that some short proofs are eliminated from the search space; see also \cite{BoolosCut}.
This includes short proofs that can be found by today's HO ATPs when the right cut/comprehension introduction steps are applied. 
The complete avoidance of cut/comprehension introduction turns certain solvable  problems into unsolvable ones; BCP is one such example. 
A hybrid approach seems to be required. 

\item Controlled cut/comprehension introduction should be considered to be a challenge for the 21st century, expanding on the progress that has been made with regards to cut-elimination since the last century. 
Machine learning may have an interesting role to play here. 

\item To some extent statement 2 is not entirely new. Earlier research on induction theorem proving, proof planning and proof methods \cite{ComputationalLogic,ScienceOfReasoning,SiekmannMelis} can be seen as pioneering work on controlled cut/comprehension introduction. 
Unfortunately this line of research did not receive the attention that it deserved at the time, due the lack of robustness and coverage of the proof planning approaches, especially when compared with the much more technically advanced FO ATPs of the time; cf.~the discussions by Bundy \cite{Critique} and Benzmüller et al.~\cite{FreshStart}. 
There is, however, successful recent related work in this area, including, e.g., work on lemma discovery for induction \cite{LemmaDiscovery}.
\end{enumerate}

The shorthand notations $ind$ and $p$, resp.~lemma $L$, introduced to solve BCP with 
HO ATPs, are actually less out of the blue than they might seem at first glance. Rather, they are indicative of a general proof method:
\begin{enumerate}
\item Introducing shorthand notation for set predicates such as $ind$ seems a good idea in general when certain inductive definitions are found in the input problem, as exemplified here by axioms A1--A3 of the Ackermann function. The formulation of such predicates is quite straightforward. 

\item The systematic introduction of related $p$-predicates also seems quite practicable. The idea is to express that results of applying an inductively defined function are always contained in the smallest $ind$-set. As noted in footnote \ref{footnote}, this sort of operation can also be viewed from an compositional (algebraic) perspective. 

\item Alternatively, these two steps can be combined and a respective analogon to lemma L can be used.
\end{enumerate}

In future work we plan to experiment with the implementation and assessment of such a general method for controlled cut/comprehension introduction in HO ATPs, in order to enable HO ATPs to find short proofs for problems like BCP.
This will include further investigation of the surprise effect revealed in this paper: By proposing appropriate shorthand notation, inspirational lemma introduction steps that would normally require further cut-introductions are now synthesized by the cut-free search procedures in HO ATPs. 
This unexpected observation deserves further attention and clarification.

The ability to find or construct short proofs plays an important role not only for FO and HO ATPs, but also for SAT and SMT solvers; see, e.g., the techniques presented by Heule et al.~\cite{HeuleKB17}. Such techniques are orthogonal to the findings reported in this paper, which exploit the gain in expressivity when moving from a less expressive logic to a more expressive one.

\section{Conclusion}
\label{Conclusion}

Holly has a strategy to never lose the contest, provided she chooses BCP as her challenge problem. 
Whatever FO ATP Folbert chooses, it will never be able to find and express a short proof for BCP in its FO calculus. 
On the flip side, if Holly chooses a state of the art HO ATP, while it might not be able to solve Folbert's challenge problem on the first days of the contest,
one day it will integrate the controlled cut/comprehension introduction techniques 
(and include further related techniques exploiting HO expressiveness), so that it 
will be able to speculate the necessary lemmas, and be able find a short proof for 
BCP (and other related challenge problems). 
The core observation is that Holly's HO ATP can steadily be improved regarding (i)~its traditional FO proof search capabilities, and (ii)~clever exploitation of its HO expressivity advantage.  
In contrast, Folbert's FO ATP is stuck with only (i), which will never lead to a solution of BCP. 
Even if the contest might (for the time being) end up in numerous draws, only Holly ever has a chance of entering heaven. 


\paragraph{Acknowledgements}{We thank the anonymous reviewers for their valuable feedback.}

\singlespacing
\bibliography{main}


\begin{thebibliography}{32}
\ifx \bisbn   \undefined \def \bisbn  #1{ISBN #1}\fi
\ifx \binits  \undefined \def \binits#1{#1}\fi
\ifx \bauthor  \undefined \def \bauthor#1{#1}\fi
\ifx \batitle  \undefined \def \batitle#1{#1}\fi
\ifx \bjtitle  \undefined \def \bjtitle#1{#1}\fi
\ifx \bvolume  \undefined \def \bvolume#1{\textbf{#1}}\fi
\ifx \byear  \undefined \def \byear#1{#1}\fi
\ifx \bissue  \undefined \def \bissue#1{#1}\fi
\ifx \bfpage  \undefined \def \bfpage#1{#1}\fi
\ifx \blpage  \undefined \def \blpage #1{#1}\fi
\ifx \burl  \undefined \def \burl#1{\textsf{#1}}\fi
\ifx \doiurl  \undefined \def \doiurl#1{\url{https://doi.org/#1}}\fi
\ifx \betal  \undefined \def \betal{\textit{et al.}}\fi
\ifx \binstitute  \undefined \def \binstitute#1{#1}\fi
\ifx \binstitutionaled  \undefined \def \binstitutionaled#1{#1}\fi
\ifx \bctitle  \undefined \def \bctitle#1{#1}\fi
\ifx \beditor  \undefined \def \beditor#1{#1}\fi
\ifx \bpublisher  \undefined \def \bpublisher#1{#1}\fi
\ifx \bbtitle  \undefined \def \bbtitle#1{#1}\fi
\ifx \bedition  \undefined \def \bedition#1{#1}\fi
\ifx \bseriesno  \undefined \def \bseriesno#1{#1}\fi
\ifx \blocation  \undefined \def \blocation#1{#1}\fi
\ifx \bsertitle  \undefined \def \bsertitle#1{#1}\fi
\ifx \bsnm \undefined \def \bsnm#1{#1}\fi
\ifx \bsuffix \undefined \def \bsuffix#1{#1}\fi
\ifx \bparticle \undefined \def \bparticle#1{#1}\fi
\ifx \barticle \undefined \def \barticle#1{#1}\fi
\bibcommenthead
\ifx \bconfdate \undefined \def \bconfdate #1{#1}\fi
\ifx \botherref \undefined \def \botherref #1{#1}\fi
\ifx \url \undefined \def \url#1{\textsf{#1}}\fi
\ifx \bchapter \undefined \def \bchapter#1{#1}\fi
\ifx \bbook \undefined \def \bbook#1{#1}\fi
\ifx \bcomment \undefined \def \bcomment#1{#1}\fi
\ifx \oauthor \undefined \def \oauthor#1{#1}\fi
\ifx \citeauthoryear \undefined \def \citeauthoryear#1{#1}\fi
\ifx \endbibitem  \undefined \def \endbibitem {}\fi
\ifx \bconflocation  \undefined \def \bconflocation#1{#1}\fi
\ifx \arxivurl  \undefined \def \arxivurl#1{\textsf{#1}}\fi
\csname PreBibitemsHook\endcsname

\bibitem{GoedelProofLength}
\begin{bchapter}
\bauthor{\bsnm{Gödel}, \binits{K.}}:
\bctitle{{Über die Länge von Beweisen}}.
In: \beditor{\bsnm{Menger}, \binits{K.}},
\beditor{\bsnm{G{\"o}del}, \binits{K.}},
\beditor{\bsnm{Wald}, \binits{A.}} (eds.)
\bbtitle{{Ergebnisse} Eines Mathematischen {Kolloquiums}: {Heft} X7:
  1934-1935},
pp. \bfpage{23}--\blpage{24}.
\bpublisher{Franz Deuticke},
\blocation{Vienna}
(\byear{1936})
\end{bchapter}
\endbibitem

\bibitem{BCI}
\begin{barticle}
\bauthor{\bsnm{Boolos}, \binits{G.}}:
\batitle{A curious inference}.
\bjtitle{J. Philos. Log.}
\bvolume{16}(\bissue{1}),
\bfpage{1}--\blpage{12}
(\byear{1987}).
\doiurl{10.1007/BF00250612}
\end{barticle}
\endbibitem

\bibitem{Eprover}
\begin{bchapter}
\bauthor{\bsnm{Schulz}, \binits{S.}},
\bauthor{\bsnm{Cruanes}, \binits{S.}},
\bauthor{\bsnm{Vukmirovic}, \binits{P.}}:
\bctitle{Faster, higher, stronger: {E} 2.3}.
In: \beditor{\bsnm{Fontaine}, \binits{P.}} (ed.)
\bbtitle{{CADE} 2019}.
\bsertitle{Lecture Notes in Computer Science},
vol. \bseriesno{11716},
pp. \bfpage{495}--\blpage{507}.
\bpublisher{Springer},
\blocation{Cham}
(\byear{2019}).
\doiurl{10.1007/978-3-030-29436-6_29}
\end{bchapter}
\endbibitem

\bibitem{CVC5}
\begin{bchapter}
\bauthor{\bsnm{Barbosa}, \binits{H.}},
\bauthor{\bsnm{Barrett}, \binits{C.W.}},
\bauthor{\bsnm{Brain}, \binits{M.}},
\bauthor{\bsnm{Kremer}, \binits{G.}},
\bauthor{\bsnm{Lachnitt}, \binits{H.}},
\bauthor{\bsnm{Mann}, \binits{M.}},
\bauthor{\bsnm{Mohamed}, \binits{A.}},
\bauthor{\bsnm{Mohamed}, \binits{M.}},
\bauthor{\bsnm{Niemetz}, \binits{A.}},
\bauthor{\bsnm{N{\"{o}}tzli}, \binits{A.}},
\bauthor{\bsnm{Ozdemir}, \binits{A.}},
\bauthor{\bsnm{Preiner}, \binits{M.}},
\bauthor{\bsnm{Reynolds}, \binits{A.}},
\bauthor{\bsnm{Sheng}, \binits{Y.}},
\bauthor{\bsnm{Tinelli}, \binits{C.}},
\bauthor{\bsnm{Zohar}, \binits{Y.}}:
\bctitle{cvc5: {A} versatile and industrial-strength {SMT} solver}.
In: \beditor{\bsnm{Fisman}, \binits{D.}},
\beditor{\bsnm{Rosu}, \binits{G.}} (eds.)
\bbtitle{Tools and Algorithms for the Construction and Analysis of Systems -
  28th International Conference, {TACAS} 2022, Held as Part of the European
  Joint Conferences on Theory and Practice of Software, {ETAPS} 2022, Munich,
  Germany, April 2-7, 2022, Proceedings, Part {I}}.
\bsertitle{Lecture Notes in Computer Science},
vol. \bseriesno{13243},
pp. \bfpage{415}--\blpage{442}.
\bpublisher{Springer},
\blocation{Cham}
(\byear{2022}).
\doiurl{10.1007/978-3-030-99524-9_24}
\end{bchapter}
\endbibitem

\bibitem{Ehoh}
\begin{barticle}
\bauthor{\bsnm{Vukmirovic}, \binits{P.}},
\bauthor{\bsnm{Blanchette}, \binits{J.}},
\bauthor{\bsnm{Cruanes}, \binits{S.}},
\bauthor{\bsnm{Schulz}, \binits{S.}}:
\batitle{Extending a brainiac prover to lambda-free higher-order logic}.
\bjtitle{International Journal on Software Tools for Technology Transfer}
\bvolume{24}(\bissue{1}),
\bfpage{67}--\blpage{87}
(\byear{2022}).
\doiurl{10.1007/s10009-021-00639-7}
\end{barticle}
\endbibitem

\bibitem{Leo3}
\begin{barticle}
\bauthor{\bsnm{Steen}, \binits{A.}},
\bauthor{\bsnm{Benzm{\"u}ller}, \binits{C.}}:
\batitle{Extensional higher-order paramodulation in {Leo-III}}.
\bjtitle{Journal of Automated Reasoning}
\bvolume{65}(\bissue{6}),
\bfpage{775}--\blpage{807}
(\byear{2021}).
\doiurl{10.1007/s10817-021-09588-x}
\end{barticle}
\endbibitem

\bibitem{Zipperpin}
\begin{barticle}
\bauthor{\bsnm{Bentkamp}, \binits{A.}},
\bauthor{\bsnm{Blanchette}, \binits{J.}},
\bauthor{\bsnm{Tourret}, \binits{S.}},
\bauthor{\bsnm{Vukmirovic}, \binits{P.}},
\bauthor{\bsnm{Waldmann}, \binits{U.}}:
\batitle{Superposition with lambdas}.
\bjtitle{Journal of Automated Reasoning}
\bvolume{65}(\bissue{7}),
\bfpage{893}--\blpage{940}
(\byear{2021}).
\doiurl{10.1007/s10817-021-09595-y}
\end{barticle}
\endbibitem

\bibitem{LostProof}
\begin{bchapter}
\bauthor{\bsnm{Benzm{\"u}ller}, \binits{C.}},
\bauthor{\bsnm{Kerber}, \binits{M.}}:
\bctitle{A lost proof}.
In: \bbtitle{Proceedings of the IJCAR 2001 Workshop: Future Directions in
  Automated Reasoning},
\bconflocation{Siena, Italy},
pp. \bfpage{13}--\blpage{24}
(\byear{2001}).
\burl{https://www.inf.ed.ac.uk/publications/online/0046/b040.pdf}
\end{bchapter}
\endbibitem

\bibitem{ackermann1928}
\begin{barticle}
\bauthor{\bsnm{Ackermann}, \binits{W.}}:
\batitle{{Zum Hilbertschen Aufbau der reellen Zahlen}}.
\bjtitle{Mathematische Annalen}
\bvolume{99}(\bissue{1}),
\bfpage{118}--\blpage{133}
(\byear{1928}).
\doiurl{10.1007/BF01459088}
\end{barticle}
\endbibitem

\bibitem{Sut22-IGPL}
\begin{barticle}
\bauthor{\bsnm{Sutcliffe}, \binits{G.}}:
\batitle{The logic languages of the {TPTP} world}.
\bjtitle{Logic Journal of the IGPL}
(\byear{2022}).
\doiurl{10.1093/jigpal/jzac068}
\end{barticle}
\endbibitem

\bibitem{THF}
\begin{barticle}
\bauthor{\bsnm{Sutcliffe}, \binits{G.}},
\bauthor{\bsnm{Benzm{\"u}ller}, \binits{C.}}:
\batitle{Automated reasoning in higher-order logic using the {TPTP THF}
  infrastructure}.
\bjtitle{Journal of Formalized Reasoning}
\bvolume{3}(\bissue{1}),
\bfpage{1}--\blpage{27}
(\byear{2010}).
\doiurl{10.6092/issn.1972-5787/1710}
\end{barticle}
\endbibitem

\bibitem{Omega}
\begin{bchapter}
\bauthor{\bsnm{Autexier}, \binits{S.}},
\bauthor{\bsnm{Benzm{\"u}ller}, \binits{C.}},
\bauthor{\bsnm{Dietrich}, \binits{D.}},
\bauthor{\bsnm{Siekmann}, \binits{J.}}:
\bctitle{{OMEGA}: Resource-adaptive processes in an automated reasoning
  systems}.
In: \beditor{\bsnm{Crocker}, \binits{M.W.}},
\beditor{\bsnm{Siekmann}, \binits{J.}} (eds.)
\bbtitle{Resource-Adaptive Cognitive Processes}.
\bsertitle{Cognitive Technologies},
pp. \bfpage{389}--\blpage{423}.
\bpublisher{Springer},
\blocation{Berlin, Heidelberg}
(\byear{2010}).
\doiurl{10.1007/978-3-540-89408-7_17}
\end{bchapter}
\endbibitem

\bibitem{Mizar}
\begin{bchapter}
\bauthor{\bsnm{Bancerek}, \binits{G.}},
\bauthor{\bsnm{Bylinski}, \binits{C.}},
\bauthor{\bsnm{Grabowski}, \binits{A.}},
\bauthor{\bsnm{Kornilowicz}, \binits{A.}},
\bauthor{\bsnm{Matuszewski}, \binits{R.}},
\bauthor{\bsnm{Naumowicz}, \binits{A.}},
\bauthor{\bsnm{Pak}, \binits{K.}},
\bauthor{\bsnm{Urban}, \binits{J.}}:
\bctitle{Mizar: State-of-the-art and beyond}.
In: \beditor{\bsnm{Kerber}, \binits{M.}},
\beditor{\bsnm{Carette}, \binits{J.}},
\beditor{\bsnm{Kaliszyk}, \binits{C.}},
\beditor{\bsnm{Rabe}, \binits{F.}},
\beditor{\bsnm{Sorge}, \binits{V.}} (eds.)
\bbtitle{Intelligent Computer Mathematics - International Conference, {CICM}
  2015, Washington, DC, USA, July 13-17, 2015, Proceedings}.
\bsertitle{Lecture Notes in Computer Science},
vol. \bseriesno{9150},
pp. \bfpage{261}--\blpage{279}.
\bpublisher{Springer},
\blocation{Cham}
(\byear{2015}).
\doiurl{10.1007/978-3-319-20615-8_17}
\end{bchapter}
\endbibitem

\bibitem{BoolosOmegaMizar}
\begin{bchapter}
\bauthor{\bsnm{Benzm{\"u}ller}, \binits{C.}},
\bauthor{\bsnm{Brown}, \binits{C.}}:
\bctitle{The curious inference of {B}oolos in {MIZAR} and {OMEGA}}.
In: \beditor{\bsnm{Matuszewski}, \binits{R.}},
\beditor{\bsnm{Zalewska}, \binits{A.}} (eds.)
\bbtitle{From Insight to Proof -- Festschrift in Honour of {Andrzej}
  {Trybulec}}.
\bsertitle{Studies in Logic, Grammar, and Rhetoric},
vol. \bseriesno{10(23)},
pp. \bfpage{299}--\blpage{388}.
\bpublisher{The University of Bialystok},
\blocation{Poland}
(\byear{2007}).
\bcomment{\href{http://mizar.org/trybulec65/20.pdf}{http://mizar.org/trybulec65/20.pdf}}
\end{bchapter}
\endbibitem

\bibitem{KetlandAFP}
\begin{botherref}
\oauthor{\bsnm{Ketland}, \binits{J.}}:
Boolos's curious inference in {Isabelle/HOL}.
Archive of Formal Proofs,
1--19
(2022).
\href{https://isa-afp.org/entries/Boolos_Curious_Inference.html}{https://isa-afp.org/entries/Boolos\_Curious\_Inference.html}
\end{botherref}
\endbibitem

\bibitem{Isabelle}
\begin{bbook}
\bauthor{\bsnm{Nipkow}, \binits{T.}},
\bauthor{\bsnm{Paulson}, \binits{L.C.}},
\bauthor{\bsnm{Wenzel}, \binits{M.}}:
\bbtitle{{Isabelle/HOL} - {A} Proof Assistant for Higher-Order Logic}.
\bsertitle{Lecture Notes in Computer Science},
vol. \bseriesno{2283}.
\bpublisher{Springer},
\blocation{Berlin, Heidelberg}
(\byear{2002}).
\doiurl{10.1007/3-540-45949-9}
\end{bbook}
\endbibitem

\bibitem{Church40}
\begin{barticle}
\bauthor{\bsnm{Church}, \binits{A.}}:
\batitle{A formulation of the simple theory of types}.
\bjtitle{Journal of Symbolic Logic}
\bvolume{5}(\bissue{2}),
\bfpage{56}--\blpage{68}
(\byear{1940}).
\doiurl{10.2307/2266170}
\end{barticle}
\endbibitem

\bibitem{Andrews}
\begin{bbook}
\bauthor{\bsnm{Andrews}, \binits{P.B.}}:
\bbtitle{An Introduction to Mathematical Logic and Type Theory}.
\bsertitle{Applied Logic Series}.
\bpublisher{Springer},
\blocation{Netherlands}
(\byear{2002}).
\doiurl{10.1007/978-94-015-9934-4}
\end{bbook}
\endbibitem

\bibitem{SEP}
\begin{bchapter}
\bauthor{\bsnm{Benzm{\"u}ller}, \binits{C.}},
\bauthor{\bsnm{Andrews}, \binits{P.}}:
\bctitle{Church's type theory}.
In: \beditor{\bsnm{Zalta}, \binits{E.N.}} (ed.)
\bbtitle{The Stanford Encyclopedia of Philosophy},
\bedition{Summer 2019} edn.,
pp. \bfpage{1}--\blpage{62}.
\bpublisher{Metaphysics Research Lab, Stanford University},
\blocation{CA, USA}
(\byear{2019}).
\bcomment{\href{https://plato.stanford.edu/entries/type-theory-church/}{https://plato.stanford.edu/entries/type-theory-church/}}
\end{bchapter}
\endbibitem

\bibitem{Andrews89}
\begin{barticle}
\bauthor{\bsnm{Andrews}, \binits{P.B.}}:
\batitle{On connections and higher-order logic}.
\bjtitle{Journal of Automated Reasoning}
\bvolume{5}(\bissue{3}),
\bfpage{257}--\blpage{291}
(\byear{1989}).
\doiurl{10.1007/BF00248320}
\end{barticle}
\endbibitem

\bibitem{Leo2}
\begin{barticle}
\bauthor{\bsnm{Benzm{\"u}ller}, \binits{C.}},
\bauthor{\bsnm{Sultana}, \binits{N.}},
\bauthor{\bsnm{Paulson}, \binits{L.C.}},
\bauthor{\bsnm{Theiss}, \binits{F.}}:
\batitle{The higher-order prover {LEO-II}}.
\bjtitle{Journal of Automated Reasoning}
\bvolume{55}(\bissue{4}),
\bfpage{389}--\blpage{404}
(\byear{2015}).
\doiurl{10.1007/s10817-015-9348-y}
\end{barticle}
\endbibitem

\bibitem{TPTP}
\begin{barticle}
\bauthor{\bsnm{Sutcliffe}, \binits{G.}}:
\batitle{The {TPTP} problem library and associated infrastructure - from {CNF}
  to {TH0}, {TPTP} v6.4.0}.
\bjtitle{Journal of Automated Reasoning}
\bvolume{59}(\bissue{4}),
\bfpage{483}--\blpage{502}
(\byear{2017}).
\doiurl{10.1007/s10817-017-9407-7}
\end{barticle}
\endbibitem

\bibitem{KV13}
\begin{bchapter}
\bauthor{\bsnm{Kov{\'{a}}cs}, \binits{L.}},
\bauthor{\bsnm{Voronkov}, \binits{A.}}:
\bctitle{First-order theorem proving and {Vampire}}.
In: \beditor{\bsnm{Sharygina}, \binits{N.}},
\beditor{\bsnm{Veith}, \binits{H.}} (eds.)
\bbtitle{Computer Aided Verification - 25th International Conference, {CAV}
  2013, Saint Petersburg, Russia, July 13-19, 2013. Proceedings}.
\bsertitle{Lecture Notes in Computer Science},
vol. \bseriesno{8044},
pp. \bfpage{1}--\blpage{35}.
\bpublisher{Springer},
\blocation{Berlin, Heidelberg}
(\byear{2013}).
\doiurl{10.1007/978-3-642-39799-8_1}
\end{bchapter}
\endbibitem

\bibitem{SST14}
\begin{bchapter}
\bauthor{\bsnm{Stump}, \binits{A.}},
\bauthor{\bsnm{Sutcliffe}, \binits{G.}},
\bauthor{\bsnm{Tinelli}, \binits{C.}}:
\bctitle{{StarExec}: {A} cross-community infrastructure for logic solving}.
In: \beditor{\bsnm{Demri}, \binits{S.}},
\beditor{\bsnm{Kapur}, \binits{D.}},
\beditor{\bsnm{Weidenbach}, \binits{C.}} (eds.)
\bbtitle{Automated Reasoning - 7th International Joint Conference, {IJCAR}
  2014, Held as Part of the Vienna Summer of Logic, {VSL} 2014, Vienna,
  Austria, July 19-22, 2014. Proceedings}.
\bsertitle{Lecture Notes in Computer Science},
vol. \bseriesno{8562},
pp. \bfpage{367}--\blpage{373}.
\bpublisher{Springer},
\blocation{Cham}
(\byear{2014}).
\doiurl{10.1007/978-3-319-08587-6_28}
\end{bchapter}
\endbibitem

\bibitem{BoolosCut}
\begin{barticle}
\bauthor{\bsnm{Boolos}, \binits{G.}}:
\batitle{Don't eliminate cut}.
\bjtitle{Journal of Philosophical Logic}
\bvolume{13}(\bissue{4}),
\bfpage{373}--\blpage{378}
(\byear{1984}).
\doiurl{10.1007/BF00247711}
\end{barticle}
\endbibitem

\bibitem{ComputationalLogic}
\begin{bbook}
\bauthor{\bsnm{Boyer}, \binits{R.S.}},
\bauthor{\bsnm{Moore}, \binits{J.S.}}:
\bbtitle{A Computational Logic}.
\bsertitle{{ACM} monograph series}.
\bpublisher{Academic Press},
\blocation{London}
(\byear{1979})
\end{bbook}
\endbibitem

\bibitem{ScienceOfReasoning}
\begin{bchapter}
\bauthor{\bsnm{Bundy}, \binits{A.}}:
\bctitle{A science of reasoning}.
In: \beditor{\bsnm{Lassez}, \binits{J.}},
\beditor{\bsnm{Plotkin}, \binits{G.D.}} (eds.)
\bbtitle{Computational Logic - Essays in Honor of Alan Robinson},
pp. \bfpage{178}--\blpage{198}.
\bpublisher{The {MIT} Press},
\blocation{Cambridge, MA, USA}
(\byear{1991}).
\bcomment{\href{https://mitpress.mit.edu/9780262121569/}{https://mitpress.mit.edu/9780262121569/}}
\end{bchapter}
\endbibitem

\bibitem{SiekmannMelis}
\begin{barticle}
\bauthor{\bsnm{Melis}, \binits{E.}},
\bauthor{\bsnm{Siekmann}, \binits{J.H.}}:
\batitle{Knowledge-based proof planning}.
\bjtitle{Artificial Intelligence}
\bvolume{115}(\bissue{1}),
\bfpage{65}--\blpage{105}
(\byear{1999}).
\doiurl{10.1016/S0004-3702(99)00076-4}
\end{barticle}
\endbibitem

\bibitem{Critique}
\begin{bchapter}
\bauthor{\bsnm{Bundy}, \binits{A.}}:
\bctitle{A critique of proof planning}.
In: \beditor{\bsnm{Kakas}, \binits{A.C.}},
\beditor{\bsnm{Sadri}, \binits{F.}} (eds.)
\bbtitle{Computational Logic: Logic Programming and Beyond, Essays in Honour of
  {Robert} {A.} {Kowalski}, Part {II}}.
\bsertitle{Lecture Notes in Computer Science},
vol. \bseriesno{2408},
pp. \bfpage{160}--\blpage{177}.
\bpublisher{Springer},
\blocation{Berlin, Heidelberg}
(\byear{2002}).
\doiurl{10.1007/3-540-45632-5_7}
\end{bchapter}
\endbibitem

\bibitem{FreshStart}
\begin{bchapter}
\bauthor{\bsnm{Benzm{\"u}ller}, \binits{C.}},
\bauthor{\bsnm{Meier}, \binits{A.}},
\bauthor{\bsnm{Melis}, \binits{E.}},
\bauthor{\bsnm{Pollet}, \binits{M.}},
\bauthor{\bsnm{Sorge}, \binits{V.}}:
\bctitle{Proof planning: A fresh start?}
In: \bbtitle{Proceedings of the {IJCAR} 2001 Workshop: Future Directions in
  Automated Reasoning},
\bconflocation{Siena, Italy},
pp. \bfpage{25}--\blpage{37}
(\byear{2001}).
\bcomment{\href{https://www.researchgate.net/publication/319393675}{https://www.researchgate.net/publication/319393675}}
\end{bchapter}
\endbibitem

\bibitem{LemmaDiscovery}
\begin{bchapter}
\bauthor{\bsnm{Johansson}, \binits{M.}}:
\bctitle{Lemma discovery for induction - {A} survey}.
In: \beditor{\bsnm{Kaliszyk}, \binits{C.}},
\beditor{\bsnm{Brady}, \binits{E.C.}},
\beditor{\bsnm{Kohlhase}, \binits{A.}},
\beditor{\bsnm{Coen}, \binits{C.S.}} (eds.)
\bbtitle{Intelligent Computer Mathematics - 12th International Conference,
  {CICM} 2019, Prague, Czech Republic, July 8-12, 2019, Proceedings}.
\bsertitle{Lecture Notes in Computer Science},
vol. \bseriesno{11617},
pp. \bfpage{125}--\blpage{139}.
\bpublisher{Springer},
\blocation{Cham}
(\byear{2019}).
\doiurl{10.1007/978-3-030-23250-4_9}
\end{bchapter}
\endbibitem

\bibitem{HeuleKB17}
\begin{bchapter}
\bauthor{\bsnm{Heule}, \binits{M.J.H.}},
\bauthor{\bsnm{Kiesl}, \binits{B.}},
\bauthor{\bsnm{Biere}, \binits{A.}}:
\bctitle{Short proofs without new variables}.
In: \beditor{\bparticle{de} \bsnm{Moura}, \binits{L.}} (ed.)
\bbtitle{Automated Deduction - {CADE} 26 - 26th International Conference on
  Automated Deduction, Gothenburg, Sweden, August 6-11, 2017, Proceedings}.
\bsertitle{Lecture Notes in Computer Science},
vol. \bseriesno{10395},
pp. \bfpage{130}--\blpage{147}.
\bpublisher{Springer},
\blocation{Cham}
(\byear{2017}).
\doiurl{10.1007/978-3-319-63046-5\_9}
\end{bchapter}
\endbibitem

\end{thebibliography}
\newpage
\appendix


\section{BCP with $ind$ and $p$ encoded in TPTP THF}
\label{sec:Boolos1.p}
Encoding of BCP with $ind$ and $p$ in TPTP THF syntax \cite{THF}; see also the~\href{https://github.com/cbenzmueller/BoolosCuriousInference-ATP/tree/main/Boolos1.p}{online} sources at \url{https://github.com/cbenzmueller/BoolosCuriousInference-ATP/tree/main/Boolos1.p}.\\

\fvset{frame=single,numbers=none,xleftmargin=.5mm,numbersep=3pt,fontsize=\relsize{-2}}
\begin{Verbatim}
% Declarations
thf(e_decl,type, e: $i ).
thf(s_decl,type, s: $i > $i ).
thf(d_decl,type, d: $i > $o ).
thf(f_decl,type, f: $i > $i > $i ).

% Auxiliary declarations
thf(ind_decl,type, ind: ( $i > $o ) > $o ).
thf(p_decl,type, p: $i > $i > $o ).

% Axioms 
thf(a1,axiom,
    ! [N: $i] :
      ( ( f @ N @ e )
      = ( s @ e ) ) ).
thf(a2,axiom,
    ! [Y: $i] :
      ( ( f @ e @ ( s @ Y ) )
      = ( s @ ( s @ ( f @ e @ Y ) ) ) ) ).
thf(a3,axiom,
    ! [X: $i,Y: $i] :
      ( ( f @ ( s @ X ) @ ( s @ Y ) )
      = ( f @ X @ ( f @ ( s @ X ) @ Y ) ) ) ).
thf(a4,axiom,
    d @ e ).
thf(a5,axiom,
    ! [X: $i] :
      ( ( d @ X )
     => ( d @ ( s @ X ) ) ) ).

% Shorthand notations  
thf(ind_def,axiom,
    ( ind
    = ( ^ [Q: $i > $o] :
          ( ( Q @ e )
          & ! [X: $i] :
              ( ( Q @ X )
             => ( Q @ ( s @ X ) ) ) ) ) ) ).
   
thf(p_def,axiom,
    ( p
    = ( ^ [X: $i,Y: $i] :
          ( ^ [X2: $i] :
            ! [X3: $i > $o] :
              ( ( ind @ X3 )
             => ( X3 @ X2 ) )
          @ ( f @ X @ Y ) ) ) ) ).

% Conjecture 
thf(conj_0,conjecture,
    d @ ( f @ ( s @ ( s @ ( s @ ( s @ e ) ) ) ) @ ( s @ ( s @ ( s @ ( s @ e ) ) ) ) ) ).
\end{Verbatim}

\pagebreak
\section{IDV proof by E for BCP with $ind$ and $p$}
\label{sec:Boolos1IDV.pdf}
\href{https://github.com/cbenzmueller/BoolosCuriousInference-ATP/tree/main/Boolos1IDV.pdf}{IDV display}  of the resolution proof by E for BCP with $ind$ and $p$ (highlighted in red are the dependencies of clause $c\_0\_52$); the resolution proof by E is presented in full detail in Appendix \ref{sec:Boolos1.proof}.

\includegraphics[height=.89\textheight]{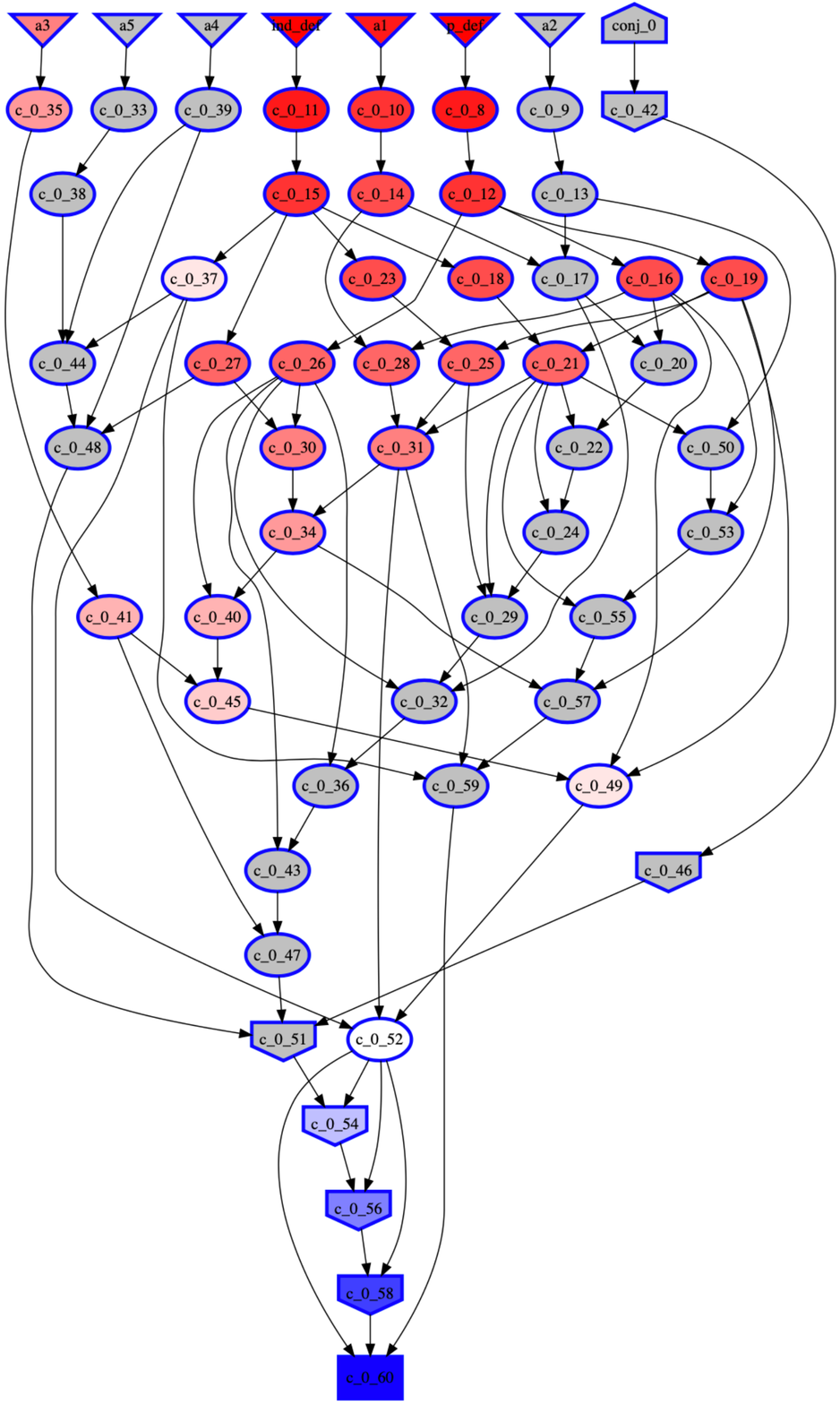}

\pagebreak
\section{Plain proof by E for BCP with $ind$ and $p$} 
\label{sec:Boolos1.proof}

See also the~\href{https://github.com/cbenzmueller/BoolosCuriousInference-ATP/tree/main/Boolos1.proof}
{online} sources at \url{https://github.com/cbenzmueller/BoolosCuriousInference-ATP/tree/main/Boolos1.proof}. 
The exact call to prover E was (where {\footnotesize \verb|/tmp/AjzF3EheOy/SOT_78srga|} refers to the input file for E as presented in Appendix \ref{sec:Boolos1.p}): 
\vskip.5em
{\footnotesize
\begin{verbatim}
 /home/tptp/Systems/E---3.0/eprover-ho --delete-bad-limit=2000000000 --definitional-cnf=24 
 -s --print-statistics -R --print-version --proof-object --auto-schedule=8 --cpu-limit=60 
 /tmp/AjzF3EheOy/SOT_78srga
\end{verbatim}
}
\vskip1em

\fvset{frame=single,numbers=none,xleftmargin=.5mm,numbersep=3pt,fontsize=\relsize{-3}}
\begin{Verbatim}
%------------------------------------------------------------------------------
% File     : E---3.0
% Problem  : SOT_78srga : v?.?.?
% Transfm  : none
% Format   : tptp:raw
% Command  : run_E %s %d THM

% Computer : quoll.cs.miami.edu
% Model    : x86_64 x86_64
% CPU      : Intel(R) Xeon(R) CPU E5-4610 0 @ 2.40GHz
% Memory   : 128720MB
% OS       : Linux 3.10.0-1160.71.1.el7.x86_64
% CPULimit : 60s
% WCLimit  : 0s
% DateTime : Fri Jul 29 06:14:27 EDT 2022

% Result   : Theorem 14.38s 1.87s
% Output   : CNFRefutation 14.38s
% Verified : 
% SZS Type : Refutation
%            Derivation depth      :   16
%            Number of leaves      :   16
% Syntax   : Number of formulae    :   69 (  26 unt;   8 typ;   0 def)
%            Number of atoms       :  141 (  12 equ;   0 cnn)
%            Maximal formula atoms :    7 (   2 avg)
%            Number of connectives :  545 (  54   ~;  61   |;   7   &; 416   @)
%                                         (   2 <=>;   5  =>;   0  <=;   0 <~>)
%            Maximal formula depth :   13 (   6 avg)
%            Number of types       :    2 (   0 usr)
%            Number of type conns  :   33 (  33   >;   0   *;   0   +;   0  <<)
%            Number of symbols     :   10 (   8 usr;   2 con; 0-3 aty)
%            Number of variables   :   77 (   4   ^  73   !;   0   ?;  77   :)

% Comments : 
%------------------------------------------------------------------------------
thf(decl_22,type,
    e: $i ).

thf(decl_23,type,
    s: $i > $i ).

thf(decl_24,type,
    d: $i > $o ).

thf(decl_25,type,
    f: $i > $i > $i ).

thf(decl_26,type,
    ind: ( $i > $o ) > $o ).

thf(decl_27,type,
    p: $i > $i > $o ).

thf(decl_28,type,
    esk1_1: ( $i > $o ) > $i ).

thf(decl_29,type,
    epred1_2: $i > $i > $i > $o ).

thf(p_def,axiom,
    ( p
    = ( ^ [X3: $i,X2: $i] :
          ( ^ [X5: $i] :
            ! [X6: $i > $o] :
              ( ( ind @ X6 )
             => ( X6 @ X5 ) )
          @ ( f @ X3 @ X2 ) ) ) ),
    file('/tmp/AjzF3EheOy/SOT_78srga',p_def) ).

thf(a2,axiom,
    ! [X2: $i] :
      ( ( f @ e @ ( s @ X2 ) )
      = ( s @ ( s @ ( f @ e @ X2 ) ) ) ),
    file('/tmp/AjzF3EheOy/SOT_78srga',a2) ).

thf(a1,axiom,
    ! [X1: $i] :
      ( ( f @ X1 @ e )
      = ( s @ e ) ),
    file('/tmp/AjzF3EheOy/SOT_78srga',a1) ).

thf(ind_def,axiom,
    ( ind
    = ( ^ [X4: $i > $o] :
          ( ( X4 @ e )
          & ! [X3: $i] :
              ( ( X4 @ X3 )
             => ( X4 @ ( s @ X3 ) ) ) ) ) ),
    file('/tmp/AjzF3EheOy/SOT_78srga',ind_def) ).

thf(a5,axiom,
    ! [X3: $i] :
      ( ( d @ X3 )
     => ( d @ ( s @ X3 ) ) ),
    file('/tmp/AjzF3EheOy/SOT_78srga',a5) ).

thf(a3,axiom,
    ! [X3: $i,X2: $i] :
      ( ( f @ ( s @ X3 ) @ ( s @ X2 ) )
      = ( f @ X3 @ ( f @ ( s @ X3 ) @ X2 ) ) ),
    file('/tmp/AjzF3EheOy/SOT_78srga',a3) ).

thf(a4,axiom,
    d @ e,
    file('/tmp/AjzF3EheOy/SOT_78srga',a4) ).

thf(conj_0,conjecture,
    d @ ( f @ ( s @ ( s @ ( s @ ( s @ e ) ) ) ) @ ( s @ ( s @ ( s @ ( s @ e ) ) ) ) ),
    file('/tmp/AjzF3EheOy/SOT_78srga',conj_0) ).

thf(c_0_8,plain,
    ! [X13: $i,X14: $i] :
      ( ( p @ X13 @ X14 )
    <=> ! [X6: $i > $o] :
          ( ( ind @ X6 )
         => ( X6 @ ( f @ X13 @ X14 ) ) ) ),
    inference(fof_simplification,[status(thm)],
      [inference(fof_simplification,[status(thm)],[p_def])]) ).

thf(c_0_9,plain,
    ! [X16: $i] :
      ( ( f @ e @ ( s @ X16 ) )
      = ( s @ ( s @ ( f @ e @ X16 ) ) ) ),
    inference(variable_rename,[status(thm)],[a2]) ).

thf(c_0_10,plain,
    ! [X15: $i] :
      ( ( f @ X15 @ e )
      = ( s @ e ) ),
    inference(variable_rename,[status(thm)],[a1]) ).

thf(c_0_11,plain,
    ! [X12: $i > $o] :
      ( ( ind @ X12 )
    <=> ( ( X12 @ e )
        & ! [X3: $i] :
            ( ( X12 @ X3 )
           => ( X12 @ ( s @ X3 ) ) ) ) ),
    inference(fof_simplification,[status(thm)],
      [inference(fof_simplification,[status(thm)],[ind_def])]) ).

thf(c_0_12,plain,
    ! [X24: $i,X25: $i,X26: $i > $o,X27: $i,X28: $i] :
      ( ( ~ ( p @ X24 @ X25 )
        | ~ ( ind @ X26 )
        | ( X26 @ ( f @ X24 @ X25 ) ) )
      & ( ( ind @ ( epred1_2 @ X27 @ X28 ) )
        | ( p @ X27 @ X28 ) )
      & ( ~ ( epred1_2 @ X27 @ X28 @ ( f @ X27 @ X28 ) )
        | ( p @ X27 @ X28 ) ) ),
    inference(distribute,[status(thm)],[inference(shift_quantors,[status(thm)],
     [inference(skolemize,[status(esa)],
      [inference(variable_rename,[status(thm)],[inference(shift_quantors,[status(thm)],
       [inference(fof_nnf,[status(thm)],[c_0_8])])])])])]) ).

thf(c_0_13,plain,
    ! [X1: $i] :
      ( ( f @ e @ ( s @ X1 ) )
      = ( s @ ( s @ ( f @ e @ X1 ) ) ) ),
    inference(split_conjunct,[status(thm)],[c_0_9]) ).

thf(c_0_14,plain,
    ! [X1: $i] :
      ( ( f @ X1 @ e )
      = ( s @ e ) ),
    inference(split_conjunct,[status(thm)],[c_0_10]) ).

thf(c_0_15,plain,
    ! [X20: $i > $o,X21: $i,X22: $i > $o] :
      ( ( ( X20 @ e )
        | ~ ( ind @ X20 ) )
      & ( ~ ( X20 @ X21 )
        | ( X20 @ ( s @ X21 ) )
        | ~ ( ind @ X20 ) )
      & ( ( X22 @ ( esk1_1 @ X22 ) )
        | ~ ( X22 @ e )
        | ( ind @ X22 ) )
      & ( ~ ( X22 @ ( s @ ( esk1_1 @ X22 ) ) )
        | ~ ( X22 @ e )
        | ( ind @ X22 ) ) ),
    inference(distribute,[status(thm)],
      [inference(shift_quantors,[status(thm)],[inference(skolemize,[status(esa)],
       [inference(variable_rename,[status(thm)],[inference(shift_quantors,[status(thm)],
        [inference(fof_nnf,[status(thm)],[c_0_11])])])])])]) ).

thf(c_0_16,plain,
    ! [X1: $i,X2: $i] :
      ( ( p @ X1 @ X2 )
      | ~ ( epred1_2 @ X1 @ X2 @ ( f @ X1 @ X2 ) ) ),
    inference(split_conjunct,[status(thm)],[c_0_12]) ).

thf(c_0_17,plain,
    ( ( f @ e @ ( s @ e ) )
    = ( s @ ( s @ ( s @ e ) ) ) ),
    inference(spm,[status(thm)],[c_0_13,c_0_14]) ).

thf(c_0_18,plain,
    ! [X1: $i,X4: $i > $o] :
      ( ( X4 @ ( s @ X1 ) )
      | ~ ( X4 @ X1 )
      | ~ ( ind @ X4 ) ),
    inference(split_conjunct,[status(thm)],[c_0_15]) ).

thf(c_0_19,plain,
    ! [X1: $i,X2: $i] :
      ( ( ind @ ( epred1_2 @ X1 @ X2 ) )
      | ( p @ X1 @ X2 ) ),
    inference(split_conjunct,[status(thm)],[c_0_12]) ).

thf(c_0_20,plain,
    ( ( p @ e @ ( s @ e ) )
    | ~ ( epred1_2 @ e @ ( s @ e ) @ ( s @ ( s @ ( s @ e ) ) ) ) ),
    inference(spm,[status(thm)],[c_0_16,c_0_17]) ).

thf(c_0_21,plain,
    ! [X1: $i,X2: $i,X3: $i] :
      ( ( epred1_2 @ X1 @ X2 @ ( s @ X3 ) )
      | ( p @ X1 @ X2 )
      | ~ ( epred1_2 @ X1 @ X2 @ X3 ) ),
    inference(spm,[status(thm)],[c_0_18,c_0_19]) ).

thf(c_0_22,plain,
    ( ( p @ e @ ( s @ e ) )
    | ~ ( epred1_2 @ e @ ( s @ e ) @ ( s @ ( s @ e ) ) ) ),
    inference(spm,[status(thm)],[c_0_20,c_0_21]) ).

thf(c_0_23,plain,
    ! [X4: $i > $o] :
      ( ( X4 @ e )
      | ~ ( ind @ X4 ) ),
    inference(split_conjunct,[status(thm)],[c_0_15]) ).

thf(c_0_24,plain,
    ( ( p @ e @ ( s @ e ) )
    | ~ ( epred1_2 @ e @ ( s @ e ) @ ( s @ e ) ) ),
    inference(spm,[status(thm)],[c_0_22,c_0_21]) ).

thf(c_0_25,plain,
    ! [X1: $i,X2: $i] :
      ( ( epred1_2 @ X1 @ X2 @ e )
      | ( p @ X1 @ X2 ) ),
    inference(spm,[status(thm)],[c_0_23,c_0_19]) ).

thf(c_0_26,plain,
    ! [X1: $i,X2: $i,X4: $i > $o] :
      ( ( X4 @ ( f @ X1 @ X2 ) )
      | ~ ( p @ X1 @ X2 )
      | ~ ( ind @ X4 ) ),
    inference(split_conjunct,[status(thm)],[c_0_12]) ).

thf(c_0_27,plain,
    ! [X4: $i > $o] :
      ( ( X4 @ ( esk1_1 @ X4 ) )
      | ( ind @ X4 )
      | ~ ( X4 @ e ) ),
    inference(split_conjunct,[status(thm)],[c_0_15]) ).

thf(c_0_28,plain,
    ! [X1: $i] :
      ( ( p @ X1 @ e )
      | ~ ( epred1_2 @ X1 @ e @ ( s @ e ) ) ),
    inference(spm,[status(thm)],[c_0_16,c_0_14]) ).

thf(c_0_29,plain,
    p @ e @ ( s @ e ),
    inference(csr,[status(thm)],[inference(spm,[status(thm)],[c_0_24,c_0_21]),c_0_25]) ).

thf(c_0_30,plain,
    ! [X1: $i,X4: $i > $o] :
      ( ( X4 @ ( f @ X1 @ ( esk1_1 @ ( p @ X1 ) ) ) )
      | ( ind @ ( p @ X1 ) )
      | ~ ( p @ X1 @ e )
      | ~ ( ind @ X4 ) ),
    inference(spm,[status(thm)],[c_0_26,c_0_27]) ).

thf(c_0_31,plain,
    ! [X1: $i] : ( p @ X1 @ e ),
    inference(csr,[status(thm)],[inference(spm,[status(thm)],[c_0_28,c_0_21]),c_0_25]) ).

thf(c_0_32,plain,
    ! [X4: $i > $o] :
      ( ( X4 @ ( s @ ( s @ ( s @ e ) ) ) )
      | ~ ( ind @ X4 ) ),
    inference(rw,[status(thm)],[inference(spm,[status(thm)],[c_0_26,c_0_29]),c_0_17]) ).

thf(c_0_33,plain,
    ! [X19: $i] :
      ( ~ ( d @ X19 )
      | ( d @ ( s @ X19 ) ) ),
    inference(variable_rename,[status(thm)],[inference(fof_nnf,[status(thm)],[a5])]) ).

thf(c_0_34,plain,
    ! [X1: $i,X4: $i > $o] :
      ( ( X4 @ ( f @ X1 @ ( esk1_1 @ ( p @ X1 ) ) ) )
      | ( ind @ ( p @ X1 ) )
      | ~ ( ind @ X4 ) ),
    inference(cn,[status(thm)],[inference(rw,[status(thm)],[c_0_30,c_0_31])]) ).

thf(c_0_35,plain,
    ! [X17: $i,X18: $i] :
      ( ( f @ ( s @ X17 ) @ ( s @ X18 ) )
      = ( f @ X17 @ ( f @ ( s @ X17 ) @ X18 ) ) ),
    inference(variable_rename,[status(thm)],[a3]) ).

thf(c_0_36,plain,
    ! [X1: $i,X4: $i > $o] :
      ( ( X4 @ ( f @ X1 @ ( s @ ( s @ ( s @ e ) ) ) ) )
      | ~ ( ind @ ( p @ X1 ) )
      | ~ ( ind @ X4 ) ),
    inference(spm,[status(thm)],[c_0_26,c_0_32]) ).

thf(c_0_37,plain,
    ! [X4: $i > $o] :
      ( ( ind @ X4 )
      | ~ ( X4 @ ( s @ ( esk1_1 @ X4 ) ) )
      | ~ ( X4 @ e ) ),
    inference(split_conjunct,[status(thm)],[c_0_15]) ).

thf(c_0_38,plain,
    ! [X1: $i] :
      ( ( d @ ( s @ X1 ) )
      | ~ ( d @ X1 ) ),
    inference(split_conjunct,[status(thm)],[c_0_33]) ).

thf(c_0_39,plain,
    d @ e,
    inference(split_conjunct,[status(thm)],[a4]) ).

thf(c_0_40,plain,
    ! [X1: $i,X2: $i,X4: $i > $o] :
      ( ( X4 @ ( f @ X1 @ ( f @ X2 @ ( esk1_1 @ ( p @ X2 ) ) ) ) )
      | ( ind @ ( p @ X2 ) )
      | ~ ( ind @ ( p @ X1 ) )
      | ~ ( ind @ X4 ) ),
    inference(spm,[status(thm)],[c_0_26,c_0_34]) ).

thf(c_0_41,plain,
    ! [X1: $i,X2: $i] :
      ( ( f @ ( s @ X1 ) @ ( s @ X2 ) )
      = ( f @ X1 @ ( f @ ( s @ X1 ) @ X2 ) ) ),
    inference(split_conjunct,[status(thm)],[c_0_35]) ).

thf(c_0_42,negated_conjecture,
    ~ ( d @ ( f @ ( s @ ( s @ ( s @ ( s @ e ) ) ) ) @ ( s @ ( s @ ( s @ ( s @ e ) ) ) ) ) ),
    inference(fof_simplification,[status(thm)],[inference(assume_negation,[status(cth)],[conj_0])]) ).

thf(c_0_43,plain,
    ! [X1: $i,X2: $i,X4: $i > $o] :
      ( ( X4 @ ( f @ X1 @ ( f @ X2 @ ( s @ ( s @ ( s @ e ) ) ) ) ) )
      | ~ ( ind @ ( p @ X2 ) )
      | ~ ( ind @ ( p @ X1 ) )
      | ~ ( ind @ X4 ) ),
    inference(spm,[status(thm)],[c_0_26,c_0_36]) ).

thf(c_0_44,plain,
    ( ( ind @ d )
    | ~ ( d @ ( esk1_1 @ d ) ) ),
    inference(cn,[status(thm)],[inference(rw,[status(thm)],
      [inference(spm,[status(thm)],[c_0_37,c_0_38]),c_0_39])]) ).

thf(c_0_45,plain,
    ! [X1: $i,X4: $i > $o] :
      ( ( X4 @ ( f @ ( s @ X1 ) @ ( s @ ( esk1_1 @ ( p @ ( s @ X1 ) ) ) ) ) )
      | ( ind @ ( p @ ( s @ X1 ) ) )
      | ~ ( ind @ ( p @ X1 ) )
      | ~ ( ind @ X4 ) ),
    inference(spm,[status(thm)],[c_0_40,c_0_41]) ).

thf(c_0_46,negated_conjecture,
    ~ ( d @ ( f @ ( s @ ( s @ ( s @ ( s @ e ) ) ) ) @ ( s @ ( s @ ( s @ ( s @ e ) ) ) ) ) ),
    inference(split_conjunct,[status(thm)],[c_0_42]) ).

thf(c_0_47,plain,
    ! [X1: $i,X4: $i > $o] :
      ( ( X4 @ ( f @ ( s @ X1 ) @ ( s @ ( s @ ( s @ ( s @ e ) ) ) ) ) )
      | ~ ( ind @ ( p @ ( s @ X1 ) ) )
      | ~ ( ind @ ( p @ X1 ) )
      | ~ ( ind @ X4 ) ),
    inference(spm,[status(thm)],[c_0_43,c_0_41]) ).

thf(c_0_48,plain,
    ind @ d,
    inference(cn,[status(thm)],[inference(rw,[status(thm)],
      [inference(spm,[status(thm)],[c_0_44,c_0_27]),c_0_39])]) ).

thf(c_0_49,plain,
    ! [X1: $i] :
      ( ( p @ ( s @ X1 ) @ ( s @ ( esk1_1 @ ( p @ ( s @ X1 ) ) ) ) )
      | ( ind @ ( p @ ( s @ X1 ) ) )
      | ~ ( ind @ ( p @ X1 ) ) ),
    inference(csr,[status(thm)],[inference(spm,[status(thm)],[c_0_16,c_0_45]),c_0_19]) ).

thf(c_0_50,plain,
    ! [X1: $i,X2: $i,X3: $i] :
      ( ( epred1_2 @ X1 @ X2 @ ( f @ e @ ( s @ X3 ) ) )
      | ( p @ X1 @ X2 )
      | ~ ( epred1_2 @ X1 @ X2 @ ( s @ ( f @ e @ X3 ) ) ) ),
    inference(spm,[status(thm)],[c_0_21,c_0_13]) ).

thf(c_0_51,negated_conjecture,
    ( ~ ( ind @ ( p @ ( s @ ( s @ ( s @ ( s @ e ) ) ) ) ) )
    | ~ ( ind @ ( p @ ( s @ ( s @ ( s @ e ) ) ) ) ) ),
    inference(cn,[status(thm)],[inference(rw,[status(thm)],
      [inference(spm,[status(thm)],[c_0_46,c_0_47]),c_0_48])]) ).

thf(c_0_52,plain,
    ! [X1: $i] :
      ( ( ind @ ( p @ ( s @ X1 ) ) )
      | ~ ( ind @ ( p @ X1 ) ) ),
    inference(cn,[status(thm)],[inference(rw,[status(thm)],
      [inference(spm,[status(thm)],[c_0_37,c_0_49]),c_0_31])]) ).

thf(c_0_53,plain,
    ! [X1: $i] :
      ( ( p @ e @ ( s @ X1 ) )
      | ~ ( epred1_2 @ e @ ( s @ X1 ) @ ( s @ ( f @ e @ X1 ) ) ) ),
    inference(spm,[status(thm)],[c_0_16,c_0_50]) ).

thf(c_0_54,negated_conjecture,
    ~ ( ind @ ( p @ ( s @ ( s @ ( s @ e ) ) ) ) ),
    inference(spm,[status(thm)],[c_0_51,c_0_52]) ).

thf(c_0_55,plain,
    ! [X1: $i] :
      ( ( p @ e @ ( s @ X1 ) )
      | ~ ( epred1_2 @ e @ ( s @ X1 ) @ ( f @ e @ X1 ) ) ),
    inference(spm,[status(thm)],[c_0_53,c_0_21]) ).

thf(c_0_56,negated_conjecture,
    ~ ( ind @ ( p @ ( s @ ( s @ e ) ) ) ),
    inference(spm,[status(thm)],[c_0_54,c_0_52]) ).

thf(c_0_57,plain,
    ( ( p @ e @ ( s @ ( esk1_1 @ ( p @ e ) ) ) )
    | ( ind @ ( p @ e ) ) ),
    inference(csr,[status(thm)],[inference(spm,[status(thm)],[c_0_55,c_0_34]),c_0_19]) ).

thf(c_0_58,negated_conjecture,
    ~ ( ind @ ( p @ ( s @ e ) ) ),
    inference(spm,[status(thm)],[c_0_56,c_0_52]) ).

thf(c_0_59,plain,
    ind @ ( p @ e ),
    inference(cn,[status(thm)],[inference(rw,[status(thm)],
      [inference(spm,[status(thm)],[c_0_37,c_0_57]),c_0_31])]) ).

thf(c_0_60,negated_conjecture,
    $false,
    inference(cn,[status(thm)],[inference(rw,[status(thm)],
      [inference(spm,[status(thm)],[c_0_58,c_0_52]),c_0_59])]),[proof] ).

%------------------------------------------------------------------------------
%----ORIGINAL SYSTEM OUTPUT
% Running higher-order theorem proving
% Running: /home/tptp/Systems/E---3.0/eprover-ho --delete-bad-limit=2000000000
%   --definitional-cnf=24 -s --print-statistics -R --print-version --proof-object
%   --auto-schedule=8 --cpu-limit=60 /tmp/AjzF3EheOy/SOT_78srga
% # Version: 3.0pre003-ho
% # partial match(1): HSMSSMSSSSSNFFN
% # Preprocessing class: HSMSSMSSSSSNHFN.
% # Scheduled 4 strats onto 8 cores with 60 seconds (480 total)
% # Starting new_ho_10 with 300s (5) cores
% # Starting new_bool_1 with 60s (1) cores
% # Starting full_lambda_5 with 60s (1) cores
% # Starting new_ho_10_unif with 60s (1) cores
% # new_ho_10 with pid 63439 completed with status 0
% # Result found by new_ho_10
% # partial match(1): HSMSSMSSSSSNFFN
% # Preprocessing class: HSMSSMSSSSSNHFN.
% # Scheduled 4 strats onto 8 cores with 60 seconds (480 total)
% # Starting new_ho_10 with 300s (5) cores
% # No SInE strategy applied
% # Search class: HGUSS-FFSF21-MHFFFSBN
% # partial match(2): HGUSS-FFSF11-MHHFFSBN
% # Scheduled 6 strats onto 5 cores with 300 seconds (300 total)
% # Starting new_ho_10_unif with 163s (1) cores
% # Starting new_ho_10 with 31s (1) cores
% # Starting lpo8_s with 28s (1) cores
% # Starting sh5l with 28s (1) cores
% # Starting sh2lt with 28s (1) cores
% # new_ho_10 with pid 63446 completed with status 0
% # Result found by new_ho_10
% # partial match(1): HSMSSMSSSSSNFFN
% # Preprocessing class: HSMSSMSSSSSNHFN.
% # Scheduled 4 strats onto 8 cores with 60 seconds (480 total)
% # Starting new_ho_10 with 300s (5) cores
% # No SInE strategy applied
% # Search class: HGUSS-FFSF21-MHFFFSBN
% # partial match(2): HGUSS-FFSF11-MHHFFSBN
% # Scheduled 6 strats onto 5 cores with 300 seconds (300 total)
% # Starting new_ho_10_unif with 163s (1) cores
% # Starting new_ho_10 with 31s (1) cores
% # Preprocessing time       : 0.004 s
% # Presaturation interreduction done
% 
% # Proof found!
% # SZS status Theorem
% # SZS output start CNFRefutation
% See solution above
% # Parsed axioms                        : 14
% # Removed by relevancy pruning/SinE    : 0
% # Initial clauses                      : 19
% # Removed in clause preprocessing      : 6
% # Initial clauses in saturation        : 13
% # Processed clauses                    : 4451
% # ...of these trivial                  : 952
% # ...subsumed                          : 2217
% # ...remaining for further processing  : 1282
% # Other redundant clauses eliminated   : 0
% # Clauses deleted for lack of memory   : 0
% # Backward-subsumed                    : 77
% # Backward-rewritten                   : 139
% # Generated clauses                    : 54638
% # ...of the previous two non-redundant : 49054
% # ...aggressively subsumed             : 0
% # Contextual simplify-reflections      : 95
% # Paramodulations                      : 54638
% # Factorizations                       : 0
% # NegExts                              : 0
% # Equation resolutions                 : 0
% # Propositional unsat checks           : 0
% #    Propositional check models        : 0
% #    Propositional check unsatisfiable : 0
% #    Propositional clauses             : 0
% #    Propositional clauses after purity: 0
% #    Propositional unsat core size     : 0
% #    Propositional preprocessing time  : 0.000
% #    Propositional encoding time       : 0.000
% #    Propositional solver time         : 0.000
% #    Success case prop preproc time    : 0.000
% #    Success case prop encoding time   : 0.000
% #    Success case prop solver time     : 0.000
% # Current number of processed clauses  : 1053
% #    Positive orientable unit clauses  : 443
% #    Positive unorientable unit clauses: 10
% #    Negative unit clauses             : 4
% #    Non-unit-clauses                  : 596
% # Current number of unprocessed clauses: 44213
% # ...number of literals in the above   : 85482
% # Current number of archived formulas  : 0
% # Current number of archived clauses   : 229
% # Clause-clause subsumption calls (NU) : 180166
% # Rec. Clause-clause subsumption calls : 178829
% # Non-unit clause-clause subsumptions  : 2351
% # Unit Clause-clause subsumption calls : 4024
% # Rewrite failures with RHS unbound    : 0
% # BW rewrite match attempts            : 9295
% # BW rewrite match successes           : 149
% # Condensation attempts                : 4451
% # Condensation successes               : 0
% # Termbank termtop insertions          : 4271812
% 
% # -------------------------------------------------
% # User time                : 1.728 s
% # System time              : 0.068 s
% # Total time               : 1.795 s
% # Maximum resident set size: 1844 pages
% 
% # -------------------------------------------------
% # User time                : 8.644 s
% # System time              : 0.424 s
% # Total time               : 9.067 s
% # Maximum resident set size: 1732 pages
% 
%------------------------------------------------------------------------------
\end{Verbatim}

\pagebreak
\section{BCP with $ind$ and $p$ as comprehension instances encoded in TPTP THF}
\label{sec:Boolos1Alt.p}

See also the~\href{https://github.com/cbenzmueller/BoolosCuriousInference-ATP/tree/main/Boolos1alt.p}
{online} sources at \url{https://github.com/cbenzmueller/BoolosCuriousInference-ATP/tree/main/Boolos1alt.p}. \\

\fvset{frame=single,numbers=none,xleftmargin=.5mm,numbersep=3pt,fontsize=\relsize{-2}}
\begin{Verbatim}
% Declarations
thf(e_decl,type, e: $i ).
thf(s_decl,type, s: $i > $i ).
thf(d_decl,type, d: $i > $o ).
thf(f_decl,type, f: $i > $i > $i ).

% Auxiliary declarations
thf(ind_decl,type, ind: ( $i > $o ) > $o ).
thf(p_decl,type, p: $i > $i > $o ).

% Axioms 
thf(a1,axiom,
    ! [N: $i] :
      ( ( f @ N @ e )
      = ( s @ e ) ) ).
thf(a2,axiom,
    ! [Y: $i] :
      ( ( f @ e @ ( s @ Y ) )
      = ( s @ ( s @ ( f @ e @ Y ) ) ) ) ).
thf(a3,axiom,
    ! [X: $i,Y: $i] :
      ( ( f @ ( s @ X ) @ ( s @ Y ) )
      = ( f @ X @ ( f @ ( s @ X ) @ Y ) ) ) ).
thf(a4,axiom,
    d @ e ).
thf(a5,axiom,
    ! [X: $i] :
      ( ( d @ X )
     => ( d @ ( s @ X ) ) ) ).

% Shorthand notations   
thf(ind_def,definition,
    ! [Q: $i > $o] :
      ( ( ind @ Q )
      = ( ( Q @ e )
        & ! [X: $i] :
            ( ( Q @ X )
           => ( Q @ ( s @ X ) ) ) ) ) ).
  
thf(p_def,definition,
    ! [X: $i,Y: $i] :
      ( ( p @ X @ Y )
      = ( ^ [Z: $i] : ( ! [Q: $i > $o] :
            ( ( ind @ Q )
           => ( Q @ Z ) ) ) @ ( f @ X @ Y ) ) ) ).

% Conjecture 
thf(conj_0,conjecture,
    d @ ( f @ ( s @ ( s @ ( s @ ( s @ e ) ) ) ) @ ( s @ ( s @ ( s @ ( s @ e ) ) ) ) ) ).
\end{Verbatim}

\section{Plain proof by Leo-III for BCP with $ind$ and $p$ as comprehension instances} 
\label{sec:Boolos1Alt.proof}

Leo-III's proof is provided~\href{https://github.com/cbenzmueller/BoolosCuriousInference-ATP/tree/main/Boolos1Alt.proof}{online} at: 
\url{https://github.com/cbenzmueller/BoolosCuriousInference-ATP/tree/main/Boolos1Alt.proof}.



\section{BCP with L encoded in TPTP THF}
\label{sec:BoolosComp.p}
Encoding of BCP with L in TPTP THF syntax; cf.~also the~\href{https://github.com/cbenzmueller/BoolosCuriousInference-ATP/tree/main/BoolosComp.p}
{online} sources: \\

\fvset{frame=single,numbers=none,xleftmargin=.5mm,numbersep=3pt,fontsize=\relsize{-2}}
\begin{Verbatim}
% Declarations
thf(e_decl,type, e: $i ).
thf(s_decl,type, s: $i > $i ).
thf(d_decl,type, d: $i > $o ).
thf(f_decl,type, f: $i > $i > $i ).

% Axioms 
thf(a1,axiom,
    ! [N: $i] :
      ( ( f @ N @ e )
      = ( s @ e ) ) ).
thf(a2,axiom,
    ! [Y: $i] :
      ( ( f @ e @ ( s @ Y ) )
      = ( s @ ( s @ ( f @ e @ Y ) ) ) ) ).
thf(a3,axiom,
    ! [X: $i,Y: $i] :
      ( ( f @ ( s @ X ) @ ( s @ Y ) )
      = ( f @ X @ ( f @ ( s @ X ) @ Y ) ) ) ).
thf(a4,axiom,
    d @ e ).
thf(a5,axiom,
    ! [X: $i] :
      ( ( d @ X )
     => ( d @ ( s @ X ) ) ) ).

% Shorthand notation p as comprehension instance 
thf(p_def,axiom,
    ? [P: $i > $i > $o] :
      ( P
      = ( ^ [X: $i,Y: $i] :
            ( ^ [Z: $i] :
              ! [R: $i > $o] :
                ( ( ^ [Q: $i > $o] :
                      ( ( Q @ e )
                      & ! [X: $i] :
                          ( ( Q @ X )
                         => ( Q @ ( s @ X ) ) ) )
                  @ R )
               => ( R @ Z ) )
            @ ( f @ X @ Y ) ) ) ) ).

% Conjecture
thf(conj_0,conjecture,
    d @ ( f @ ( s @ ( s @ ( s @ ( s @ e ) ) ) ) @ ( s @ ( s @ ( s @ ( s @ e ) ) ) ) ) ).
\end{Verbatim}

\section{Plain proof by E for BCP with L} 
\label{sec:BoolosComp.proof}

E's proof is provided~\href{https://github.com/cbenzmueller/BoolosCuriousInference-ATP/tree/main/Boolos1Alt.proof}{online} at: 
\url{https://github.com/cbenzmueller/BoolosCuriousInference-ATP/tree/main/BoolosComp.proof}.



\end{document}